\documentclass{amsart}

\usepackage{amssymb}
\usepackage{amscd}
\usepackage{amsmath}
\usepackage{latexsym}
\usepackage{verbatim}
\usepackage[latin1]{inputenc}
\usepackage[dvips]{graphics}

\usepackage[all]{xy}
\CompileMatrices

\newdir{ >}{{}*!/-10pt/\dir{>}}

{\begin{list}{}{%
\settowidth{\labelwidth}{\textsf{#1}}%
\setlength{\leftmargin}{\labelwidth}
\addtolength{\leftmargin}{0pt}}}%
{\end{list}}

\newtheorem{Thm}{Theorem}[section]

\newtheorem{Prop}[Thm]{Proposition}

\theoremstyle{definition}

\theoremstyle{remark}
\newtheorem{Rem}[Thm]{Remark}

\newcommand{\Prf}{\noindent\textit{Proof. }}
\newcommand{\PrfOf}[1]{\noindent\textit{Proof of #1.}}

\newcommand{\bbN}{\mathbb{N}}
\newcommand{\bbZ}{\mathbb{Z}}
\newcommand{\bbQ}{\mathbb{Q}}
\newcommand{\bbR}{\mathbb{R}}
\newcommand{\bbC}{\mathbb{C}}

\newcommand{\nspace}{\hspace*{-.1em}}
\newcommand{\nnspace}{\hspace*{-.05em}}

\newcommand{\sign}{\mathop{\rm sign}\nolimits}
\newcommand{\signxG}[2]{\sign_{\nnspace #1}^{\nnspace{\scriptscriptstyle #2}}}

\newcommand{\id}{\mathop{\rm i\hspace*{-.03em}d}\nolimits}
\newcommand{\incl}{\mathop{\rm incl}\nolimits}

\newcommand{\On}{{\rm O}}
\newcommand{\SO}{{\rm SO}}
\newcommand{\SU}{{\rm SU}}
\newcommand{\SL}{{\rm SL}}

\newcommand{\ie}{{\sl i.e.\ }}

\newcommand{\oursetminus}{\nspace\smallsetminus\nnspace}

\newcommand{\point}{\!\cdot\!}

\newcommand{\topo}{{\rm t}\nnspace{\rm o}\nnspace{\rm p}}

\newcommand{\onto}{\,-\!\!\!\!\twoheadrightarrow}

\newcommand{\EG}{\hspace*{.07em}\underline{\hspace*{-.07em}EG\hspace*{-.23em}}\hspace*{.23em}}

\newcommand{\soverline}[1]{\hspace*{.2em}\overline{\hspace*{-.2em}#1\hspace*{-.1em}}\hspace*{.1em}}

\begin{document}


\title[Homotopy invariance of higher signatures and $3$-manifold groups]{Homotopy
invariance of higher signatures\\ and $3$-manifold groups}

\author{Michel MATTHEY, Herv\'e OYONO-OYONO and Wolfgang PITSCH}

\address{University of Lausanne, IGAT, Bâtiment BCH, EPFL, CH-1015 Lausanne, Switzerland}

\urladdr{http://igat.epfl.ch/matthey/}

\email{michel.matthey@unil.ch}

\address{University Blaise Pascal in Clermont-Ferrand, Departement of Mathematics,
Complexe Scientifique des C\'ezeaux, F-63177 Aubi\`ere Cedex, France}

\email{oyono@math.univ-bpclermont.fr}

\address{University of Geneva, Mathematics Section, 2--4 rue du Li\`evre, Case postale 240,
CH-1211 Genève 24, Switzerland}

\email{wolfgang.pitsch@math.unige.ch}



\date{December 12, 2004}



\begin{abstract}
For closed oriented manifolds, we establish oriented homotopy invariance of higher signatures that come
from the fundamental group of a large class of orientable $3$-manifolds, including the ``piecewise
geometric'' ones in the sense of Thurston. In particular, this class, that will be carefully described,
is the class of all orientable $3$-manifolds if the Thurston Geometrization Conjecture is true. In fact,
for this type of groups, we show that the Baum-Connes Conjecture With Coefficients holds. The non-oriented
case is also discussed.
\end{abstract}


\maketitle


\section{Introduction and statement of the main results}
\label{s-Intro}

\medbreak


We assume all manifolds to be non-empty, pointed (\ie we fix a base-point), second countable,
Hausdorff and smooth. Given a closed connected oriented manifold $M^{m}$ of dimension $m$, let $[M]$
denote either orientation classes in $H_{m}(M;\bbQ)$ and in $H^{m}(M;\bbZ)$, and let $\mathcal{L}_
{\!M}\in H^{4*}(M;\bbQ)$ be the Hirzebruch $L$-class of $M$, which is defined as a suitable rational
polynomial in the Pontrjagin classes of $M$ (see~\cite[pp.~11--12]{Hirz} or~\cite[Ex.~III.11.15]{LawMich}).
Denote the usual Kronecker pairing for $M$, with rational coefficients, by
$$
\left<\,.\,,.\,\right>\colon H^{*}(M;\bbQ)\times H_{*}(M;\bbQ)\longrightarrow\bbQ\,.
$$
If $M$ is of dimension $m=4k$, then the Hirzebruch Signature Theorem (see~\cite[Thm.~8.2.2]{Hirz}
or~\cite[p.~133]{LawMich}) says that the rational number $\left<\mathcal{L}_{\!M},[M]\right>$ is
the signature of the cup product quadratic form
$$
H^{2k}(M;\bbZ)\otimes H^{2k}(M;\bbZ)\longrightarrow H^{4k}(M;\bbZ)=\bbZ\point[M]\cong\bbZ\,,\quad
(x,y)\longmapsto x\cup y\,.
$$
As a consequence, $\left<\mathcal{L}_{\!M},[M]\right>$ is an oriented homotopy invariant of $M$
(among closed connected oriented manifolds, hence of the same dimension $4k$). In 1965, Sergei Petrovich Novikov
proposed the following conjecture, now known as the \emph{Novikov Conjecture} or as the \emph{Novikov
Higher Signature Conjecture}\,: Let $G$ be a discrete group, let $BG$ be its classifying space, and
let $\alpha\in H^{*}(BG;\bbQ)\cong H^{*}(G;\bbQ)$ be a prescribed rational cohomology class of $BG$.
Now, for a closed connected oriented manifold $M^{m}$ (with $m$ arbitrary) and for a continuous map
$f\colon M\longrightarrow BG$, consider the \emph{$\alpha$-higher signature (coming from $G$)}
$$
\signxG{\alpha}{G}(M,f):=\big<f^{*}(\alpha)\cup\mathcal{L}_{\!M},[M]\big>\;\in\;\bbQ\,,
$$
where $f^{*}\colon H^{*}(BG;\bbQ)\longrightarrow H^{*}(M;\bbQ)$ is induced by $f$. Then, the
conjecture predicts that the rational number $\signxG{\alpha}{G}(M,f)$ is an oriented homotopy
invariant of the pair $(M,f)$, in the precise sense that $\signxG{\alpha}{G}(N,g)=\signxG{\alpha}
{G}(M,f)$ whenever $N^{n}$ is a second closed connected oriented manifold equipped with a continuous
map $g\colon N\longrightarrow BG$, and such that there exists a homotopy equivalence $h\colon M
\stackrel{\simeq}{\longrightarrow}N$ preserving the orientation, that is, $h_{*}[M]=[N]$ in $H_{m}
(N;\bbQ)$ (automatically, $m=n$), and with $g\circ h\simeq f$, \ie the diagram
$$
\xymatrix@C=2em@R=.4em{
M \ar[dd]_{h} \ar[rrd]^{f} & & \\
& \hspace*{-2.2em}\circlearrowleft\hspace*{-.94em}\raisebox{.1em}{${\scriptscriptstyle\simeq}$} & BG \\
N \ar[rru]_{g} & &
}
$$
commutes up to homotopy, as indicated. If, for a given group $G$, this holds for every rational
cohomology class $\alpha\in H^{*}(BG;\bbQ)$, then one says that $G$ verifies the Novikov Conjecture.
Of particular interest are the ``\emph{self higher signatures}'' of a closed connected oriented
manifold $M$, namely those corresponding to the case $G:=\pi_{1}(M)$, for some chosen cohomology class
$\alpha\in H^{*}(BG;\bbQ)$, with, as map $f\colon M\longrightarrow BG$, `the' classifying map of the
universal covering space $\widetilde{M}$ of $M$ (up to homotopy). Special attention is deserved by
the case where $M$ is aspherical, in which case one can take $M$ as a model for $BG$, and $f:=\id_{M}$.

\medskip

Now, fix a discrete group $G$ (countable, say). Let $K_{*}(-)$ denote complex topological $K$-homology,
with compact supports, for spaces, and let $C^{*}G$ be the maximal $C^{*}$-algebra of $G$ (a suitable
$C^{*}$-completion of the complex group algebra $\bbC G$ of $G$), whose analytical $K$-theory is denoted
by $K_{*}^{\topo}(C^{*}G)$. In \cite{Mis}, Mi\v{s}\v{c}enko defines a group homomorphism
$$
\tilde{\nu}_{*}^{G}\colon K_{*}(BG)\longrightarrow K_{*}^{\topo}(C^{*}G)
$$
and shows that if $\tilde{\nu}_{*}^{G}$ is rationally injective, \ie injective after tensoring with
$\bbQ$, then the Novikov Conjecture holds for $G$. Now, letting $C^{*}_{r}G$ be the reduced $C^{*}$-algebra
of $G$ (another suitable $C^{*}$-completion of $\bbC G$) and $\lambda^{G}\colon C^{*}G\onto C^{*}_{r}G$
the canonical surjective $*$-homomorphism, the composite
$$
\nu_{*}^{G}\colon K_{*}(BG)\stackrel{\tilde{\nu}_{*}^{G}}{\longrightarrow}K_{*}^{\topo}(C^{*}G)\stackrel
{\lambda^{G}_{*}}{\longrightarrow}K_{*}^{\topo}(C^{*}_{r}G)
$$
is called the Novikov assembly map. The so-called \emph{Strong Novikov Conjecture} for~$G$ is the statement that
$\nu_{*}^{G}$ is rationally injective, and this, again, implies the usual Novikov Conjecture. Next,
we explain the connection with the Baum-Connes Conjecture. Let $\EG$ denote the universal example for
proper actions of $G$ (in other words, up to $G$-homotopy, the classifying space for the family of finite
subgroups of~$G$); by definition, this is a locally compact Hausdorff proper (left, say) $G$-space such
that for any locally compact Hausdorff $G$-space $X$, there exists a $G$-map from $X$ to $\EG$, and any
two $G$-maps from $X$ to $\EG$ are $G$-homotopic. For instance, the universal covering $EG:=\widetilde{BG}$
of $BG$ is a model for $\EG$ when $G$ is torsion-free; the point $pt$ is a model for $\EG$ when $G$
is finite; if $G$ is a discrete subgroup of an almost connected Lie group $\Gamma$ with maximal compact
subgroup $K$, then $\Gamma/K$ is a model for $\EG$. Suppose further given a separable $G$-$C^{*}$-algebra
$\mathcal{A}$. Then, there is a suitable $G$-equivariant $K$-homology group $K_{*}^{G}(\EG;\mathcal{A})$
and a specific group homomorphism, called the Baum-Connes assembly map with coefficients in~$\mathcal{A}$,
$$
\mu_{*}^{G,\mathcal{A}}\colon K_{*}^{G}(\EG;\mathcal{A})\longrightarrow K_{*}^{\topo}(\mathcal{A}\rtimes_{r}G)\,,
$$
where $\mathcal{A}\rtimes_{r}G$ is the reduced $C^{*}$-crossed product of $\mathcal{A}$ by $G$. The group
$G$ is said to satisfy the \emph{Baum-Connes Conjecture With Coefficients} if the assembly map $\mu_{*}^{G,
\mathcal{A}}$ is an isomorphism for any separable $G$-$C^{*}$-algebra $\mathcal{A}$. If this is at least
known to be fulfilled for the $C^{*}$-algebra $\bbC$ with trivial $G$-action, then one says that $G$
verifies the \emph{Baum-Connes Conjecture} (\ie without mentioning coefficients). In this special case
where $\mathcal{A}=\bbC$ with trivial $G$-action, one has $\mathcal{A}\rtimes_{r}G=C^{*}_{r}G$ and
$K_{*}^{G}(\EG;\mathcal{A})=K_{*}^{G}(\EG)$, the $G$-equivariant $K$-homology group with $G$-compact
supports of $\EG$, and the corresponding assembly map boils down to a map
$$
\mu_{*}^{G}:=\mu_{*}^{G,\bbC}\colon K_{*}^{G}(\EG)\longrightarrow K_{*}^{\topo}(C^{*}_{r}G)\,.
$$
This is linked with the Novikov Conjecture as follows. First, since $G$ acts properly and freely
on $EG$, and since $BG\simeq G\backslash EG$, there is a canonical isomorphism
$$
K_{*}(BG)\cong K_{*}^{G}(EG)\,.
$$
Secondly, since tautologically any proper and free $G$-action is proper, there is a $G$-map
$EG\longrightarrow\EG$, unique up to $G$-homotopy, and the induced map
$$
K_{*}^{G}(EG)\longrightarrow K_{*}^{G}(\EG)
$$
is known to be rationally injective. Thirdly, the Novikov assembly map $\nu_{*}^{G}$ coincides with
the composite map
$$
K_{*}(BG)\cong K_{*}^{G}(EG)\longrightarrow K_{*}^{G}(\EG)\stackrel{\mu_{*}^{G}}{\longrightarrow}
K_{*}^{\topo}(C^{*}_{r}G)\,.
$$
It follows that if the group $G$ satisfies the Baum-Connes Conjecture (in particular, if $G$ verifies
the Baum-Connes Conjecture With Coefficients), then the Strong Novikov Conjecture holds for $G$, and
hence also the original Novikov Conjecture on higher signatures. As general references for the
Baum-Connes Conjecture and related topics, let us mention \cite{BC,BCH,MisVal,Val}.

\medskip

In this paper, we observe that so much is known about the structure of $3$-manifolds and
that the Baum-Connes Conjecture With Coefficients has been proved for such a large class
of groups, that this enables to establish the Baum-Connes Conjecture With Coefficients for
the fundamental group of any compact orientable $3$-manifold ``with a piecewise geometric
structure'', more precisely to which the famous Thurston Geometrization Conjecture applies,
namely\,:

\begin{Thm}
\label{thm-main-1}
Suppose that the Thurston Hyperbolization Conjecture is true, as for example if the Thurston Geometrization
Conjecture holds. Let $G$ be the fundamental group of an orientable $3$-manifold, compact or not, with or
without boundary. Then, the Baum-Connes Conjecture With Coefficients holds for $G$. In particular, the group
$G$ satisfies the Novikov Conjecture, \ie higher signatures coming from $G$ are oriented homotopy invariants
for closed connected oriented manifolds of arbitrary dimension.
\end{Thm}

\begin{Rem}
In Section~\ref{s-Proofs}, more details will be given about the Thurston Geometrization Conjecture and
the Thurston Hyperbolization Conjecture (see Remark~\ref{rem-Thurston-Conj} below).
\end{Rem}

\begin{Rem}
By recent outstanding results of Perelman, one might expect to have, in a near future, a complete proof of
the Thurston Geometrization Conjecture, and hence of the Thurston Hyperbolization Conjecture.
\end{Rem}

In fact, in the compact case, we have a more precise result, independently of the Thurston Hyperbolization
Conjecture\,:

\begin{Thm}
\label{thm-main-2}
Let $G$ be the fundamental group of a compact orientable $3$-manifold $M$ (possibly with boundary),
and consider a two-stage decomposition of the capped-off manifold $\widehat{M}$ of $M$, firstly,
into Kneser's prime decomposition, secondly, for each occurring closed irreducible piece with infinite
fundamental group, a Jaco-Shalen-Johannson torus decomposition. Now, consider only those pieces obtained
after the second stage and which are closed, non-Haken, non-Seifert, non-hyperbolizable and whose fundamental
group is infinite. Suppose that the fundamental groups of these very pieces all satisfy the Baum-Connes
Conjecture with Coefficients. Then, $G$ verifies the Baum-Connes Conjecture with Coefficients and the
Novikov Conjecture.
\end{Thm}

\begin{Rem}
Let $M$ be a compact $3$-manifold. The \emph{capped-off} manifold $\widehat{M}$ of $M$ is
obtained from $M$ by capping off with a compact $3$-ball each boundary component of $M$ that
is diffeomorphic to a $2$-sphere, getting this way a compact $3$-manifold $\widehat{M}$,
see~\cite[p.~25]{Hemp}. Note that $\widehat{M}$ is orientable whenever $M$ is orientable,
and that the inclusion $M\hookrightarrow\widehat{M}$ induces an isomorphism on the level
of fundamental groups.
\end{Rem}

\begin{Rem}
In Section~\ref{s-Proofs}, we will explain Kneser's and Jaco-Shalen-Johannson's decompositions.
We will also define the notions of prime, of irreducible, of Haken, of Seifert, and of hyperbolizable
$3$-manifolds.
\end{Rem}

\begin{Rem}
In particular, all ``self higher signatures'' are oriented homotopy invariants for closed connected
oriented $3$-manifolds to which Theorems \ref{thm-main-1} and \ref{thm-main-2} apply. At this point,
it is worth mentioning that all irreducible compact connected orientable $3$-manifolds with infinite
fundamental group are aspherical, as follows from the Sphere Theorem, see~\cite[p.~483]{Sco} and
\cite[Thm.~4.3]{Hemp}.
\end{Rem}

In the non-orientable compact case, we have the following result.

\begin{Thm}
\label{thm-non-or}
Let $M$ be a compact non-orientable $3$-manifold, and let $G$ be its fundamental group. Let
$M_{1},\ldots,M_{p}$ be the irreducible pieces in Kneser's (normal) prime decomposition.
Suppose, for each $i=1,\ldots,p$, that one of the following properties is fulfilled\,: either
$M_{i}$ is orientable and satisfies the hypotheses of Theorem~\ref{thm-main-2} (as for
example if Thurston Hyperbolization Conjecture is true); or $\pi_{1}(M_{i})$ is infinite
cyclic; or $M_{i}$ is non-orientable and without $2$-torsion in its fundamental group. Then,
the group $G$ satisfies the Baum-Connes Conjecture With Coefficients and the Novikov Conjecture.
\end{Thm}

\begin{Rem}
In Section~\ref{s-Proofs}, we will explain when a Kneser prime decomposition is called normal
(a property guaranteeing its uniqueness).
\end{Rem}

\begin{Rem}
The Baum-Connes Conjecture With Coefficients, hence the Novikov Conjecture, is known for the fundamental
group of any manifold of dimension $\leq 2$. So, what is done here, is extending this result up to dimension
$3$ in the orientable case, modulo the Thurston Hyperbolization Conjecture. Since, for each $n\geq 4$, every
finitely presentable group is isomorphic to the fundamental group of some closed connected orientable (smooth)
$n$-manifold (see for instance~\cite{Dehn,Mar} or~\cite{JohnWal}), a further extension one dimension up should
certainly be incomparably more difficult and seems to be, by far, out of scope at the time of writing. At this point,
we mention that by an unpublished result of Connes, Gromov and Moscovici (see however~\cite{Gro}), for closed
connected oriented manifolds of arbitrary dimension, all higher signatures coming from a discrete group $G$
and corresponding to a cohomology class lying in the subring of $H^{*}(BG;\bbQ)$ generated by the classes
of degree $\leq 2$ are oriented homotopy invariants; a complete proof is now available in~\cite[Cor.~0.3]{Math}.
\end{Rem}

\begin{Rem}
In Theorems \ref{thm-main-1}, \ref{thm-main-2} and \ref{thm-non-or}, one does not need to suppose that the
considered $3$-manifolds are smooth manifolds, but merely topological manifolds. Indeed, as is well-known,
any (second countable Hausdorff) topological manifold of dimension $\leq 3$ admits a smooth structure, which
is furthermore unique.
\end{Rem}

\begin{Rem}
If it would be known that any countable discrete group $G$ sitting in a short exact sequence
of groups
$$
1\longrightarrow H\longrightarrow G\longrightarrow\bbZ/2\longrightarrow 1\,,
$$
with $H$ satisfying the Baum-Connes Conjecture With Coefficients, verifies itself the
Baum-Connes Conjecture With Coefficients, then one could drop the condition ``orientable''
in Theorems \ref{thm-main-1} and \ref{thm-main-2} (one could also drop the first occurring
assumption of orientability in Theorem~\ref{thm-main-3} below). Indeed, suppose this is known.
Then, noticing that in both theorems there is no restriction in assuming connectedness
of the considered $3$-manifold $M$ (which is compact for~\ref{thm-main-2}), in case $M$
is non-orientable, the theorem in question applies to the orientation covering $\soverline{M}$
of $M$, which is a regular double covering of $M$ (and is itself compact for~\ref{thm-main-2}),
for which one has the fibre sequence $S^{0}\rightarrow\soverline{M}\rightarrow M$ and
therefore a short exact sequence of groups
$$
1\longrightarrow\pi_{1}(\soverline{M})\longrightarrow\pi_{1}(M)\longrightarrow\bbZ/2\longrightarrow 1\,.
$$
\end{Rem}

Before we state a consequence of our main results, recall that for a torsion-free discrete group $G$,
the \emph{Kaplansky/Idempotent Conjecture} (resp.\ the \emph{Kadison-Kaplansky Conjecture}) states
that the algebra $\bbC G$ (resp.\ $C^{*}_{r}G$) contains no non-trivial idempotent, \ie any of its
element $\varepsilon$ satisfying $\varepsilon=\varepsilon^{2}$ is equal to $0$ or $1$.

\begin{Thm}
\label{thm-main-3}
Suppose that the Thurston Hyperbolization Conjecture is true. Then, Kaplansky's Idempotent Conjecture and
the Kadison-Kaplansky Conjecture hold for any torsion-free fundamental group of an orientable $3$-manifold,
as for example for the fundamental group of any compact orientable $3$-manifold whose prime factors in
Kneser's prime decomposition all have an infinite fundamental group.
\end{Thm}

\begin{Rem}
Of course, there is a analogous statement to Theorem~\ref{thm-main-3} for all fundamental groups
to which Theorem~\ref{thm-main-2} applies, provided they are torsion-free.
\end{Rem}


\section{The proofs}
\label{s-Proofs}


We give here the proofs of Theorems \ref{thm-main-1}, \ref{thm-main-2}, \ref{thm-non-or}
and \ref{thm-main-3}.

\medskip

Before we start the proofs, we present a recollection of standard results from the
topology and geometry of $3$-manifolds. As general references on the subject, let
us cite \cite{Hemp,Sco}, and also \cite{And,Bon,Kap,Thu}.

\medskip

A $3$-manifold $M$ is called \emph{prime} if it admits no non-trivial connected sum decomposition,
\ie if $M\approx M'\#M''$, then at least one of $M'$ and $M''$ is diffeomorphic to~$S^{3}$; $M$ is
said to be \emph{irreducible} (in the sense of Hempel~\cite[p.~28]{Hemp}) if every embedded $2$-sphere
in $M$ bounds an embedded compact $3$-ball. By \cite[Lem.~3.13]{Hemp} a prime $3$-manifold is either
an $S^{2}$-bundle over $S^{1}$, or irreducible. Given an $S^{2}$-bundle $E$ over $S^{1}$, the homotopy
exact sequence of the fiber sequence $S^{2}\rightarrow E\rightarrow S^{1}$ yields that $\pi_{1}(E)$ is
infinite cyclic; if $E$ is orientable, then it is diffeomorphic to $S^{1}\times S^{2}$.

\medskip

To begin our discussion of the two-stage decomposition, we let $M$ be a compact connected $3$-manifold
(but not necessarily \emph{closed}, \ie the boundary $\partial M$ may be non-empty). By the
\emph{Kneser Prime Decomposition Theorem} (see~\cite{Kne,Mil}, or \cite[Thm.~3.15]{Hemp} where the
closeness and the orientability of $M$ are avoided, see pp.~24\,\&\,32 therein), one can decompose
$M$ as a finite connected sum of compact connected $3$-manifolds, say
$$
M\approx M_{1}\#M_{2}\#\ldots\#M_{q}\,,
$$
with each $M_{i}$ prime; we can (and will) further suppose that the decomposition is \emph{normal} in
the sense of~\cite[p.~34]{Hemp}, \ie some $M_{i}$ is diffeomorphic to $S^{1}\times S^{2}$ if and only
if $M$ is orientable. In this case, the decomposition is unique (up to reordering and diffeomorphism),
and, under the extra assumption that $M$ is orientable, each $M_{i}$ is orientable as well, see~\cite[Thm.~3.21]{Hemp}
(see also~\cite{Mil} for the orientable case). Of course, by the van Kampen Theorem, the fundamental group
of $M$ decomposes as a finite free product
$$
\pi_{1}(M)\cong\pi_{1}(M_{1})*\pi_{1}(M_{2})*\ldots*\pi_{1}(M_{q})\,.
$$
Recall that each $M_{i}$ is either an $S^{2}$-bundle over $S^{1}$, or irreducible.

\medskip

Now, we let $M$ be a compact connected $3$-manifold. In the sequel, by a \emph{surface}~$\Sigma$, we
mean a compact connected $2$-dimensional manifold (with possibly non-empty boundary $\partial
\Sigma$). Consider a surface $\Sigma$ that is either properly embedded in $M$, \ie $\partial\Sigma=
\Sigma\cap\partial M$ (transverse intersection), or embedded in $\partial M$; in case $\Sigma\subseteq
\partial M$ (so that $\Sigma$ is closed), note that `sliding' $\Sigma$ along a small collar neighbourhood
inside $M$, which is a trivial half-line bundle, we get an isotopic properly embedded surface in $M$.
The surface $\Sigma$ is called \emph{$2$-sided} if it is embedded in $\partial M$, or if it admits a
tubular neighbourhood in $M$ which is a trivial line bundle. The surface $\Sigma$ is said to be
\emph{incompressible} inside $M$ if it is $2$-sided, not diffeomorphic to the $2$-sphere nor
to a disk, and if it is \emph{$\pi_{1}$-injective}, in the sense that the inclusion $\Sigma
\hookrightarrow M$ induces a monomorphism $\pi_{1}(\Sigma)\hookrightarrow\pi_{1}(M)$. A $3$-manifold
$M$ is called \emph{$P^{2}$-irreducible} if it is irreducible and if it contains no embedded
$2$-sided real projective plane.

\medskip

A compact connected $3$-manifold $M$ is called \emph{Haken} if it is $P^{2}$-irreducible and
contains a properly embedded $2$-sided incompressible surface ($M$ is supposed to be orientable,
this amounts to require $M$ to be irreducible and to contain a properly embedded incompressible
orientable surface). By \cite[Lem.~6.7\,(i)]{Hemp}, if the compact connected $3$-manifold $M$ is orientable and
if $\partial M$ is non-empty and does not only consist of a collection of $2$-spheres, then the group
$H_{1}(M;\bbZ)$ is infinite, and in this case, \cite[Lem.~6.6]{Hemp} shows that $M$ is Haken provided
it is irreducible (the surface $F$ constructed in the proof therein indeed is orientable). A compact
connected $3$-manifold $M$ is called \emph{torus-irreducible} (or \emph{geometrically atoroidal}) if
every incompressible $2$-torus in $M$ is isotopic to a boundary component of $M$.

\medskip

For the general definition, that we will not need, of a \emph{Seifert $3$-manifold}, we refer
to \cite[pp.~428\,\&\,429]{Sco}; what we will however need is the following characterization
due to Epstein~\cite{Eps} in the compact case\,: a compact $3$-manifold $M$ is Seifert if it
admits a foliation by circles. By~\cite[Thm.~9.2]{JanNeu} (see also \cite[Thm.~1.38]{Kap}),
a deep result, a prime compact $3$-manifold $M$ with infinite fundamental group $\pi_{1}(M)$
is Seifert if and only if $\pi_{1}(M)$ contains an infinite cyclic normal subgroup, in which
case, there exists a short exact sequence of groups
$$
1\longrightarrow\bbZ\longrightarrow\pi_{1}(M)\stackrel{p}{\longrightarrow}\Gamma\longrightarrow 1\,,
$$
with $\Gamma$ standing for a discrete subgroup of the isometry group of either $S^{2}$
(the `round' $2$-sphere), of $\bbR^{2}$ (the flat Euclidean plane), or of $\mathcal{H}^{2}$
(the hyperbolic plane). This means that $\Gamma$ is a discrete subgroup of one of the following
three Lie groups (each having with exactly two connected components)\,:
$$
\On(3)\,,\qquad\bbR^{2}\rtimes\On(2)\qquad\mbox{and}\qquad\SO(2,1)\,.
$$
It will be important for us to note that for any finite subgroup $H$ of $\Gamma$, its
pre-image $p^{-1}(H)$ in $\pi_{1}(M)$ sits in a short exact sequence
$$
1\longrightarrow\bbZ\longrightarrow p^{-1}(H)\longrightarrow H\longrightarrow 1\,,
$$
and is therefore \emph{virtually cyclic}, in the sense that it contains a cyclic subgroup
(here, infinite) of finite index.

\medskip

Next, we include a short algebraic incursion. A \emph{graph of groups} $\mathcal{G}$ is a non-empty
graph $G_{\mathcal{G}}=(E_{\mathcal{G}},V_{\mathcal{G}})$ (possibly with loops, \ie with edges incident
to only one vertex, and simple, \ie with at most one loop per vertex and at most one edge joining
two distinct vertices) equipped with two families $\{G'_{e}\}_{e\in E_{\mathcal{G}}}$ and $\{G_{v}\}_
{v\in V_{\mathcal{G}}}$ of groups parameterized by the edge set $E_{\mathcal{G}}$ and the vertex set
$V_{\mathcal{G}}$, respectively, and a family $\{\iota_{e,v}\colon G'_{e}\hookrightarrow G_{v}\,|\,v
\in e\}_{e\in E_{\mathcal{G}}}$ of injective group homomorphisms, one for each pair $(e,v)\in E_
{\mathcal{G}}\times V_{\mathcal{G}}$ consisting of an edge and an adjacent vertex; the groups in
$\{G'_{e}\}_{e\in E_{\mathcal{G}}}$ and in $\{G_{v}\}_{v\in V_{\mathcal{G}}}$ are called the
\emph{edge-groups} and the \emph{vertex-groups} of $\mathcal{G}$, respectively. If the graph
of groups $\mathcal{G}$ is finite and connected (\ie if $G_{\mathcal{G}}$ is a finite connected
graph), its fundamental group $\pi_{1}(\mathcal{G})$ is a group defined, up to isomorphism, by
a finite induction process mixing the groups $G_{v}$ and $G'_{e}$, using the incidence relation
of $G_{\mathcal{G}}$ and the maps $\iota_{e,v}$, via amalgamated free products and HNN-extensions
(see~\cite[Section~5]{Serre} for details). This group $\pi_{1}(\mathcal{G})$ acts simplicially
on the graph $G_{\mathcal{G}}$, with, up to isomorphism, vertex-stabilizers $\{G_{v}\}_{v\in
V_{\mathcal{G}}}$ and edge-stabilizers $\{G'_{e}\}_{e\in E_{\mathcal{G}}}$.

\medskip

After Kneser's decomposition (or ``sphere decomposition''), there is a second decomposition that
we will need, namely the so-called \emph{JSJ-decomposition} (or ``torus decomposition''), named
after Jaco-Shalen~\cite{JacSha} and Johannson~\cite{Johan}. So, we let $M$ be an irreducible closed
connected orientable $3$-manifold. Then, there is a  minimal finite family $\{T_{\!j}\}_{j\in J}$
(possibly empty) of embedded disjoint incompressible $2$-sided closed $2$-tori that separate $M$
into a finite set $\{M_{k}\}_{k\in K}$ of irreducible compact connected orientable $3$-manifolds,
each of which is either Seifert or torus-irreducible, possibly both. (Such a family is, up to isotopy
inside $M$, unique; the finite index-sets $J$ and $K$ verify $|K|=|J|+1$.) Let us now describe
the fundamental group of $M$ using a graph of groups. It turns out that there is a graph of groups
$\mathcal{G}=\mathcal{G}_{M}$ with $E_{\mathcal{G}}=J$ and $V_{\mathcal{G}}=K$, and, for $j\in J$
and $k\in K$, $G'_{j}=\pi_{1}(T_{\!j})\cong\bbZ^{2}$, $G_{k}=\pi_{1}(M_{k})$ and $\iota_{j,k}=
\pi_{1}\big(\incl\colon T_{\!j}\hookrightarrow M_{k}\big)$, and with the incidence relation dictated
by the combinatorial configuration of the separating family of tori; moreover (and most importantly),
there is an isomorphism $\pi_{1}(M)\cong\pi_{1}(\mathcal{G})$. Indeed, this last property follows
inductively from the van Kampen Theorem.

\medskip

We also recall that an $n$-manifold $M$, possibly with non-empty boundary, is called \emph{hyperbolizable}
if its geometric interior $M\oursetminus\partial M$ admits a complete Riemannian metric for which
the sectional curvature is constant with value $-1$. In this case, $\pi_{1}(M)\cong\pi_{1}(M\oursetminus
\partial M)$ is isomorphic to a discrete subgroup of the Lie group $\SO(n,1)$ (and not necessarily
of its identity component $\SO(n,1)_{\nnspace\circ}$).

\begin{Rem}
\label{rem-Thurston-Conj}
Suppose given a closed connected orientable $3$-manifold $M$, and apply to it the following
two-stage decomposition (without necessity of first capping $M$ off). First perform Kneser's
prime decomposition; this produces finitely many pieces which are either $S^{1}\times S^{2}$
or closed irreducible manifolds. To each of the latter ones apply the JSJ-decomposition.
The \emph{Thurston Geometrization Conjecture} is the statement that the final pieces all have a
(necessary unique) geometric structure among a list of eight possible ones (in a precise and
specific sense, see~\cite{Sco,Thu}). It might well happen that one has no decomposition
to perform, for instance if one starts with~$S^{3}$. The Thurston Geometrization Conjecture
is known in all but two cases\,:
\begin{itemize}
  \item [(a)] for closed irreducible manifolds with finite fundamental group; this special case
  is known as the \emph{Thurston Elliptization Conjecture} (which is equivalent to the combination
  of the Poincaré Conjecture and of the Spherical Space Form Conjecture);

  \item [(b)] for closed, irreducible, non-Haken and non-Seifert manifolds with infinite
  fundamental group; in this case the manifold should be hyperbolizable\,: this is the content
  of the \emph{Thurston Hyperbolization Conjecture}.
\end{itemize}
There is also a more general version of the Thurston Geometrization Conjecture (that we will
not need and which is more technical to state), namely for connected orientable $3$-manifolds
that are compact (indeed, not necessarily closed). It is now known to hold in all cases, except
for the very same two `closed' cases (a) and (b).
\end{Rem}

For the proof of Theorem~\ref{thm-main-1}, we will also need the following result.

\begin{Prop}
\label{prop-non-compact}
Let $M$ be a $3$-manifold. Then, there exists a family $\{M_{n}\}_{n\in\bbN}$ of compact connected $3$-manifolds
and a family $\{f_{n}\colon M_{n}\rightarrow M\}_{n\in\bbN}$ of smooth immersions, such that each immersion
$f_{n}$ induces an injective group homomorphism $\pi_{1}(M_{n})\hookrightarrow \pi_{1}(M)$, and such that
the fundamental group of $M$ is the union of (the images of) the fundamental groups of the members
of the family, \ie
$$
\pi_{1}(M)=\bigcup_{n\in\bbN}\pi_{1}(M_{n})\,.
$$
Moreover, if $M$ is orientable, then one can further require the $M_{n}$'s to be orientable.
\end{Prop}

\Prf
First, the group $\pi_{1}(M)$ being countable, let $(g_{n})_{n\in\bbN}$ be a countable sequence
of elements of $\pi_{1}(M)$ (possibly with repetitions) such that the set $\{g_{n}\}_{n\in\bbN}$
generates $\pi_{1}(M)$. For each $n\in\bbN$, let $G_{n}:=\left<g_{1},\ldots,g_{n}\right>$ be the
subgroup of $\pi_{1}(M)$ generated by $g_{1},\ldots,g_{n}$. Fix $n\in\bbN$. Since $G_{n}$ is
finitely generated, by~\cite[Thm.~8.2]{Hemp}, it is even finitely presented. Therefore,
applying~\cite[Thm.~8.1]{Hemp} (a result due to Jaco~\cite{Jac}), we can find a compact
connected $3$-manifold $M_{n}$ and an immersion $f_{n}\colon M_{n}\rightarrow M$ such that
$(f_{n})_{*}\colon\pi_{1}(M_{n})\hookrightarrow\pi_{1}(M)$ is injective, as indicated, with
image $G_{n}$ (note that one can indeed suppose each $M_{n}$ connected). The equality $\pi_{1}
(M)=\bigcup_{n\in\bbN}\pi_{1}(M_{n})$ is now obvious. Finally, for each $n$, $M_{n}$ being
of the same dimension as $M$, and an immersion being a local homeomorphism, \cite[Ex.~3 of VIII.2.22]{Dold}
applies to~$f_{n}$ to show orientability of $M_{n}$ in case $M$ itself is orientable (note
that \cite[Prop.~VIII.2.19]{Dold} allows to incorporate successfully the case where $M_{n}$
and/or $M$ have a boundary).
\qed

\medskip

Finally, we are in position to pass to the proofs of our theorems (in disorder).

\medskip

\PrfOf{Theorem~\ref{thm-main-2}}
Clearly, for the proofs, we can suppose that the compact orientable $3$-manifold $M$ we consider
is connected, and that $M$ is capped-off, \ie that $M=\widehat{M}$. Let $G$ be the fundamental
group of $M$. From the Kneser Prime Decomposition Theorem, we have deduced a finite free product
decomposition
$$
G\cong\pi_{1}(M_{1})*\pi_{1}(M_{2})*\ldots*\pi_{1}(M_{q})\,.
$$
Since the Baum-Connes Conjecture With Coefficients is stable under forming finite free
products (see~\cite{Oyono1a,Oyono1b}), if each $\pi_{1}(M_{i})$ verifies this conjecture, then
the same holds for $G$. Since $\pi_{1}(S^{1}\times S^{2})$ is infinite cyclic, and
since the Baum-Connes Conjecture With Coefficients holds for the group $\bbZ$ (in fact,
for any countable amenable group, including all abelian groups, see~\cite{HigKas1,HigKas2}),
we can now suppose further that $M$ is irreducible. As we have explained, if $M=\widehat{M}$
is not closed, \ie if $\partial M\neq\varnothing$, then $M$ is Haken. In this case, by \cite{Oyono1b},
or \cite{BBV}, or \cite{Tu}, its fundamental group satisfies the Baum-Connes Conjecture With
Coefficients (the proof is based on the fact that a Haken manifold admits a so-called \emph{hierarchy}
in the sense of~\cite[p.~140]{Hemp} and on the results on graphs of groups we have recalled earlier).
So, we are reduced to the case where $M$ is an irreducible closed connected orientable $3$-manifold.

\smallskip

Now, we apply to $M$ a JSJ-decomposition. Earlier, in such a situation, $\pi_{1}(M)$ has been
expressed using a certain graph of groups. By~\cite{Oyono1b} again, the Baum-Connes Conjecture
With Coefficients (and also the plain Baum-Connes Conjecture, see~\cite[Thm.~5.13 in Part~I]{MisVal})
is stable under taking finite connected graphs of groups, \ie if a finite connected graph of groups
$\mathcal{G}$ has all its edge-groups $\{G'_{e}\}_{e\in E_{\mathcal{G}}}$ and vertex-groups
$\{G_{v}\}_{v\in V_{\mathcal{G}}}$ satisfying the Baum-Connes Conjecture (resp.\ With Coefficients),
then so does its fundamental group $\pi_{1}(\mathcal{G})$. As, in our case, the edge-groups are
isomorphic to the abelian group $\bbZ^{2}$, the Baum-Connes Conjecture With Coefficients holds
for them. So, it remains to deal with the vertex-groups. These are fundamental groups of compact
connected $3$-manifolds, each of which is either Seifert or torus-irreducible, possibly both.
We distinguish three cases for each of these pieces, that we call, say, $N$.

\smallskip

(1) If $N$ is Seifert, then, as we have seen, $\pi_{1}(N)$ sits in a short exact sequence
$$
1\longrightarrow\bbZ\longrightarrow\pi_{1}(N)\stackrel{p}{\longrightarrow}\Gamma\longrightarrow 1\,,
$$
with $\Gamma$ a discrete subgroup of one of the Lie groups $\On(3)$, $\bbR^{2}\rtimes\On(2)$
and $\SO(2,1)$, which are almost connected, \ie they have finitely many connected components
(here, exactly $2$). Now, consider the following facts concerning $\Gamma$\,:
\begin{itemize}
    \item [(i)] If $\Gamma$ is a discrete subgroup of the compact group $\On(3)$, then $\Gamma$ is
    finite and thus satisfies the Baum-Connes Conjecture With Coefficients (see~\cite{Julg}).
    \item [(ii)] The so-called Kasparov $\gamma$-element is equal to one for both Lie groups $\SO(2,1)$
    and $\bbR^{2}\rtimes\On(2)$. Since any discrete subgroup of an almost connected Lie group with $\gamma=1$ satisfies
    the Baum-Connes Conjecture With Coefficients, so does $\Gamma$. Given $n\geq 2$, for $\SO(n,1)$, the equality
    $\gamma=1$ is established in~\cite{Kas1}, and for $\bbR^{n}\rtimes\On(n)$, the $\gamma$-element, being
    invariant under group retractions (see~\cite{Kas1}), is the image of the $\gamma$-element of $\On(n)$,
    which, by a computation carried out in~\cite{AtS}, is equal to one as well. It could also be said that
    if $\Gamma$ is a discrete subgroup of $\SO(2,1)$ or of $\bbR^{2}\rtimes\On(2)$, then $\Gamma$ has the
    Haagerup property (see \cite[Thm. ~4.0.1 \& Prop.~6.1.5]{CCJJV} for $\Gamma\subset\SO(2,1)$,
    and~\cite{BeChVa} for $\Gamma\subset\bbR^{2}\rtimes\On(2)$, in which case $\Gamma$ is amenable)
    and then conclude by \cite{HigKas1,HigKas2}.
\end{itemize}
We have also seen that for any finite subgroup $H$ of $\Gamma$, the pre-image $p^{-1}(H)$ inside $\pi_{1}(N)$
is virtually cyclic and therefore amenable (since the class of discrete amenable groups contains abelian
groups and finite groups, and is stable under taking group extensions). By~\cite{HigKas1,HigKas2} again,
each $p^{-1}(H)$ satisfies the Baum-Connes Conjecture With Coefficients; by~\cite{Oyono2}, this is
enough to guarantee that $\pi_{1}(N)$ itself satisfies this conjecture. This is it for case~(1).

\smallskip

(2) If $N$ has finite fundamental group (hence $N$ is non-Seifert and, in fact, torus-irreducible), then
the Baum-Connes Conjecture With Coefficients is known for the finite group $\pi_{1}(N)$, as we have already
said (see~\cite{Julg}).

\smallskip

(3) If $N$ is non-Seifert with infinite fundamental group (and $N$ is then torus-irreducible), then, we
distinguish four non mutually excluding sub-cases.
\begin{itemize}
    \item [(i)] If $N$ is Haken, then, by \cite{Oyono1b}, or \cite{BBV}, or \cite{Tu}, its fundamental
    group satisfies the Baum-Connes Conjecture With Coefficients.
    \item [(ii)] If $N$ is hyperbolizable, then, as recalled earlier, $\pi_{1}(N)$ is a discrete subgroup
    of $\SO(3,1)$. As seen in (1)\,(ii) above, such a discrete subgroup satisfies the Baum-Connes Conjecture
    With Coefficients.
    \item [(iii)] If $N$ (which is non-Seifert and has infinite fundamental group) is neither Haken, nor
    hyperbolizable, then our technical hypothesis in the statement of the theorem precisely guarantees
    that $\pi_{1}(N)$ also satisfies this conjecture.
\end{itemize}
This completes our discussion of case~(3).

\smallskip

We conclude, for each considered piece $N$ obtained after the JSJ-decomposition, that, in any of these
three events (1)--(3), the group $\pi_{1}(N)$ satisfies the Baum-Connes Conjecture With Coefficients,
and consequently that so does $\pi_{1}(M)$.
\qed

\medskip

\PrfOf{Theorem~\ref{thm-main-1}}
By~\cite[Thm.~1.1]{BMP}, if a countable discrete group $G$ is the union $G=\bigcup_{n\in\bbN}G_{n}$
of a collection of subgroups all satisfying the Baum-Connes Conjecture With Coefficients, then so
does $G$. Since the fundamental group of a compact manifold is countable (at most), combining this
with Proposition~\ref{prop-non-compact}, the result follows directly from Theorem~\ref{thm-main-2};
indeed, as we have recalled, the Thurston Geometrization Conjecture implies the Thurston Hyperbolization
Conjecture, which precisely predicts that each piece obtained exactly after the second stage of the
two-stage decomposition of the statement and which is non-Seifert, non-Haken and has infinite
fundamental group is hyperbolizable.
\qed

\medskip

\PrfOf{Theorem~\ref{thm-non-or}}
We may suppose that $M$ is connected and capped-off, so that $M=\widehat{M}$. Using Kneser's
(normal) prime decomposition, we can write $M$ as
$$
M\approx M_{1}\#\ldots\#M_{p}\#M_{p+1}\#\ldots\#M_{q}
$$
with $M_{1},\ldots,M_{q}$ denoting prime compact connected $3$-manifolds (possibly non-orientable),
where $M_{1},\ldots,M_{p}$ are irreducible and $M_{p+1},\ldots,M_{q}$ are prime but not irreducible.
Therefore, $M_{p+1},\ldots,M_{q}$ are $S^{2}$-bundles over $S^{1}$ and have consequently an infinite
cyclic fundamental group, and hence verifying the Baum-Connes Conjecture With Coefficients. Now, fix
$i\in\{1,\ldots,p\}$. By assumption, either $M_{i}$ is orientable and Theorem~\ref{thm-main-2} applies
to it to show that it satisfies the Baum-Connes Conjecture With Coefficients, or $N:=M_{i}$ is an irreducible,
non-orientable, compact, connected and capped-off $3$-manifold having either infinite cyclic fundamental
group, or having no $2$-torsion in its fundamental group and with each component of $\partial M$ incompressible
in $M$ (possibly with $\partial M=\varnothing$). Let us now deal with $N$. If $\pi_{1}(N)\cong\bbZ$ then,
once again, $N$ satisfies the Baum-Connes Conjecture. So, we suppose that $\pi_{1}(N)$ is $2$-torsion-free,
but not infinite cyclic. By Kneser's Conjecture on free products, proved for instance in~\cite[Thm.~7.1]{Hemp},
since $N$ is irreducible, its fundamental group $\pi_{1}(N)$ is indecomposable with respect to free products.
This property, together with the fact that $\pi_{1}(N)$ is not infinite cyclic and does not contain $2$-torsion,
implies that \cite[Lem.~10.1]{Hemp} applies to $N$, which is capped-off. The conclusion of this result is that
$N$ is $P^{2}$-irreducible (in the notation of \cite[Lem.~10.1]{Hemp}, since $N$ is irreducible and non-orientable,
we can take $N$ as $\mathcal{P}(N)$ and the occurring homotopy sphere is diffeomorphic to $S^{3}$).
Combining~\cite[Lem.~6.7\,(ii)\,\&\,Lem.~6.6]{Hemp} for the $P^{2}$-irreducible manifold $N$, we obtain,
inside $N$, a properly embedded, $2$-sided incompressible surface $\Sigma$, which is \emph{non-separating}.
(In particular, $N$ is Haken.) Therefore, cutting $N$ along $\Sigma$, we get a compact \emph{connected}
$P^{2}$-irreducible manifold $N'$ with non-empty boundary. Invoking~\cite[Thm.~13.3]{Hemp}, we obtain a
hierarchy for $N'$ (see details in~\cite[p.~140]{Hemp}). Consequently, the argument given in~\cite{Oyono1b}
proves that the group $\pi_{1}(N')$ satisfies the Baum-Connes Conjecture With Coefficients. Now, there is
an isomorphism $\pi_{1}(N)\cong\pi_{1}(N')*_{\pi_{1}(\Sigma)}$, \ie $\pi_{1}(N)$ is an HNN-extension with
base $\pi_{1}(N')$ and over the surface group $\pi_{1}(\Sigma)$. Fundamental groups of closed surfaces
(orientable or not) are one-relator groups, so that, by~\cite{Oyono1b}, they verify the Baum-Connes Conjecture
With Coefficients. By~\cite{Oyono1b} once again, this conjecture is stable under forming HNN-extensions,
so that the conjecture holds for $\pi_{1}(N)$ too. In total, we see that each ``free factor'' in the initial
decomposition
$$
\pi_{1}(M)\cong\pi_{1}(M_{1})*\ldots*\pi_{1}(M_{p})*\pi_{1}(M_{p+1})*\ldots*\pi_{1}(M_{q})
$$
satisfies the conjecture, hence also their finite free product $\pi_{1}(M)$, still by~\cite{Oyono1b}.
\qed

\medskip

\PrfOf{Theorem~\ref{thm-main-3}}
It is standard that surjectivity of the Baum-Connes assembly map (in degree $0$) for a
torsion-free discrete group $G$ implies the Kadison-Kaplansky Conjecture for $G$,
and hence Kaplansky's Idempotent Conjecture for $G$ since $\bbC G$ is a sub-algebra
of $C^{*}_{r}G$ (see for instance~\cite[Lem.~7.2 in Part~I]{MisVal} or~\cite[Section~5]{Pusch}
for a proof). So, the first part of Theorem~\ref{thm-main-3} follows directly from
Theorem~\ref{thm-main-1}. For the second part, suppose that $G=\pi_{1}(M)$, where $M$
is a connected orientable $3$-manifold decomposed as
$$
M\approx M_{1}\#M_{2}\#\ldots\#M_{q}\,,
$$
with each $M_{i}$ a compact connected orientable prime $3$-manifold with, by assumption, infinite
fundamental group. By \cite[Thm.~9.8]{Hemp} (see also p.~170 therein), each fundamental group
$\pi_{1}(M_{i})$ is torsion-free, hence also the finite free-product $G\cong\pi_{1}(M_{1})*
\pi_{1}(M_{2})*\ldots*\pi_{1}(M_{q})$. Consequently, the first part of the theorem applies
to~$G$.
\qed




\begin{thebibliography}{00}

\bibitem{And}
M.\,T.~Anderson.
\newblock Scalar curvature and geometrization conjectures for $3$-manifolds,
\newblock In {Grove, Karsten (ed.) et al., Comparison geometry}, Cambridge University,
Math. Sci. Res. Inst. Publ. {\bf 30} (1997), 49--82.

\bibitem{AtS}
M.\,F.~Atiyah and I.\,M.~Singer.
\newblock The index of elliptic operators I,
\newblock {\em Ann. Math.} {\bf 87} (1968), 484--530.

\bibitem{BC}
P.~Baum and A.~Connes.
\newblock Chern character for discrete groups,
\newblock {\em Enseign. Math.} {\bf 46} (2000), 3--42.

\bibitem{BCH}
P.~Baum, A.~Connes, and N.~Higson.
\newblock Classifying spaces for proper actions and {$K$}-theory of group
  {$C^{*}$}-algebras,
\newblock In {\em {$C^{*}$}-algebras 1943-1993: a fifty year celebration},
  volume 167, pages 241--291, Contemporary Mathematics, 1994.

\bibitem{BMP}
P.~Baum, S.~Millington and R.~Plymen.
\newblock Local-global principle for the {B}aum-{C}onnes conjecture with coefficients,
\newblock {\em $K$-theory} {\bf 28} (2003) 1--18.

\bibitem{BBV}
C.~B\'{e}guin, H.~Bettaieb and A.~Valette.
\newblock {$K$}-theory for {$C^{*}$}-algebras of one-relator groups,
\newblock {\em $K$-theory} {\bf 16} (1999), No.3, 277--298.

\bibitem{BeChVa}
M.\,E.\,B.~Bekka, P.-A.~Cherix and A.~Valette,
\newblock Proper affine isometric actions of amenable groups,
\newblock In \emph{ Novikov conjectures, index theorems and rigidity, Vol. 2}
(Oberwolfach, 1993), 1--4,London Math. Soc. Lecture Note Ser., 227, Cambridge
Univ. Press, Cambridge, 1995.

\bibitem{Bon}
F.~Bonahon.
\newblock Geometric structures on $3$-manifolds,
\newblock In {Daverman, R. J. (ed.) et al., Handbook of geometric topology}, Elsevier
(2002), 93--164.


\bibitem{CCJJV}
P.-A.~Ch\'{e}rix, M.~Cowling, P.~Jolissaint, P.~Julg and A.~Valette.
\newblock {\em Groups with the {H}aagerup property, {G}romov's a-{T}-menability}.
\newblock Progress in Mathematics 197, Birh\"{a}user, 2001.


\bibitem{Dehn}
M.~Dehn.
\newblock \"{U}ber die {T}opologie des dreidimensionalen {R}aumes,
\newblock {\em Math. Ann.} {\bf 69} (1910), 137--168.

\bibitem{Dold}
A.~Dold.
\newblock {\em Lectures on algebraic topology}.
\newblock Classics in Mathematics, Springer Verlag, 1995.

\bibitem{Eps}
D.\,B.\,A.~Epstein.
\newblock Periodic flows on {$3$}-manifolds,
\newblock {\em Ann. of Math.} {\bf 95} (1972), 66--82.

\bibitem{Gro}
M.~Gromov.
\newblock Geometric reflections on the {N}ovikov conjecture,
\newblock In ``Novikov conjectures, index theorems and rigidity,''
  London Mathematical Society, Lecture Notes Series 226 (1995), 164--173.


\bibitem{Hemp}
J.~Hempel.
\newblock {\em $3$-manifolds}.
\newblock Annals of Math Studies 86, Princeton University Press, 1976.

\bibitem{HigKas1}
N.~Higson and G.\,G.~Kasparov.
\newblock Operator {$K$}-theory for groups which act properly and isometrically on Hilbert
space,
\newblock {\em Electron. Res. Announc. Amer. Math. Soc.} {\bf 3} (1997), 131--142 (electronic).

\bibitem{HigKas2}
N.~Higson and G.\,G.~Kasparov.
\newblock {$E$}-theory and {$K\!K$}-theory for groups which act properly and isometrically
on Hilbert space,
\newblock {\em Invent. Math.} {\bf 144} (2001), 23--74.

\bibitem{Hirz}
F.~Hirzebruch.
\newblock {\em Topological methods in algebraic geometry}.
\newblock Springer-Verlag Classics in Mathematics, Third Edition, 1995.

\bibitem{Jac}
W.~Jaco.
\newblock Finitely presented subgroups of three-manifold groups,
\newblock {\em Invent. Math.} {\bf 13} (1971), 335--346.

\bibitem{JacSha}
W.~Jaco and P.~Shalen.
\newblock Seifert fibered spaces in $3$-manifolds,
\newblock {\em Mem. Amer. Math. Soc.} {\bf 21} (1979), no. 220, 192~pp.

\bibitem{JanNeu}
M.~Jankins and W.\,D.~Neumann.
\newblock {\em Lectures on Seifert manifolds}.
\newblock Brandeis Lecture Notes 2, Brandeis University, 1983.


\bibitem{Johan}
K.~Johannson.
\newblock Homotopy equivalence of $3$-manifolds with boundaries,
\newblock {\em Springer Lecture Notes in Math.} {\bf 761}, 1979.

\bibitem{JohnWal}
F.\,E.\,A.~Johnson and J.\,P.~Walton.
\newblock Parallelizable manifolds and the fundamental group,
\newblock {\em Mathematika} {\bf 47} (2000), no. 1-2, 165--172 (2002).

\bibitem{Julg}
P.~Julg.
\newblock {$K$}-th\'eorie \'equivariante et produits crois\'es,
\newblock {\em C. R. Acad. Sci. Paris} {\bf 292} (1981), 629--632.

\bibitem{JulKas}
P.~Julg and G.\,G.~Kasparov.
\newblock L'anneau {$K\!K_{G}(\bbC,\bbC)$} pour $G=\SU(n,1)$,
\newblock {\em J. Reine Angew. Math.} {\bf 463} (1995), 99--152.

\bibitem{Kap}
M.~Kapovich.
\newblock {\em Hyperbolic manifolds and discrete groups}.
\newblock Progress in Mathematics 183, Birkh\"{a}user, 2001.

\bibitem{Kas1}
G.\,G.~Kasparov.
\newblock Lorentz groups: {$K$}-theory of unitary representations and crossed products,
\newblock {\em Soviet Math. Dokl.} {\bf 29} (1984), 256--260.

\bibitem{Kas2}
G.\,G.~Kasparov.
\newblock Equivariant {$K\!K$}-theory and the Novikov conjecture,
\newblock {\em Invent. Math.} {\bf 91} (1988), 147--201.

\bibitem{Kne}
H.~Kneser.
\newblock Geschlossene {F}l\"{a}chen in dreidimensionalen {M}annigfaltigkeiten,
\newblock {\em Jahresbericht der Deut. Math. Verein.} {\bf 38} (1929), 248--260.

\bibitem{LawMich}
H.\,B.~Lawson and M.-L.~Michelsohn.
\newblock {\em Spin geometry}.
\newblock Princeton University Press, 1989.


\bibitem{Mar}
A.~Markov.
\newblock The insolubility of the problem of homeomorphy (Russian),
\newblock {\em Dokl. Akad. Nauk SSSR} {\bf 121} (1958), 218--220.

\bibitem{Math}
V.~Mathai.
\newblock The {N}ovikov conjecture for low degree cohomology classes,
\newblock {\em Geom.\ Dedicata} {\bf 99} (2003) 1--15.

\bibitem{Mil}
J.\,W.~Milnor.
\newblock A unique factorization theorem for $3$-manifolds,
\newblock {\em Amer. J. Math.} {\bf 84} (1962), 1--7.


\bibitem{Mis}
A.\,S.~Mi\v{s}\v{c}enko.
\newblock Infinite-dimensional representations of discrete groups, and higher signatures (Russian),
\newblock {\em Izv. Akad. Nauk SSSR Ser. Mat.} {\bf 38} (1974), 81--106.

\bibitem{MisVal}
G.~Mislin and A.~Valette.
\newblock {\em Proper group actions and the {B}aum-{C}onnes conjecture}.
\newblock Course given at the Centre de Recerca Matem\`atica de Barcelona
in 2001, Advanced Course in Mathematics, CRM Barcelona, Birkh\"auser,
2003.

\bibitem{Oyono1a}
H.~Oyono-Oyono.
\newblock {B}aum-{C}onnes conjecture and group actions on trees,
\newblock {\em $K$-theory} {\bf 24} (2001), no. 2, 115--134.

\bibitem{Oyono1b}
H.~Oyono-Oyono.
\newblock La conjecture de {B}aum-{C}onnes pour les groupes agissant sur
les arbres,
\newblock {\em C. R. Acad. Sci. Paris} {\bf 326} (1998), 799--804.

\bibitem{Oyono2}
H.~Oyono-Oyono.
\newblock {B}aum-{C}onnes conjecture and extensions,
\newblock {\em J. Reine Angew. Math.} {\bf 532} (2001), 133--149.

\bibitem{Pusch}
M.~Puschnigg.
\newblock The {K}adison-{K}aplansky conjecture for word-hyperbolic groups,
\newblock {\em Invent. Math.} {\bf 149}, no. 1, 153--194, 2002.

\bibitem{Sco}
P.~Scott.
\newblock The geometries of $3$-manifolds,
\newblock {\em Bull. London Math. Soc.} {\bf 15} (1983), 401--487.

\bibitem{Serre}
J.-P.~Serre.
\newblock {\em Arbres, amalgames, $\SL_{2}$}.
\newblock Ast\'{e}risque 46, Soci\'{e}t\'{e} Math\'{e}matique de France, 1983.

\bibitem{Thu}
W.\,P.~Thurston.
\newblock {\em Three-dimensional geometry and topology, Vol.~1}.
\newblock Princeton Mathematical Series 35, Princeton University Press, 1997.

\bibitem{Tu}
J.\,L.~Tu.
\newblock The {B}aum-{C}onnes conjecture and discrete group actions on trees,
\newblock {\em $K$-theory} {\bf 17} (1999), 303--318.

\bibitem{Val}
A.~Valette.
\newblock {\em Introduction to the {B}aum-{C}onnes conjecture}.
\newblock Course given at the ETH, Z\"urich (1999), Lecture
  Notes in Mathematics ETH, Z\"urich, Birkh\"auser, 2002.

\end{thebibliography}
\end{document}